\documentclass[12pt,letterpaper]{amsart}
\usepackage{euler,latexsym, amssymb, amscd, amsfonts, xypic, url, color, epsfig, hyperref, enumerate, tikz-cd,dsfont}
\input xy
\xyoption{all}

\usepackage[T1]{fontenc}
 
 \newlength{\baseunit}               
 \newcount{\numlines}                
\newcommand{\bpf}{\noindent {\em Proof.  }}
\newcommand{\epf}{\qed \vspace{+10pt}}

\newtheorem{tm}[subsection]{Theorem}
\newtheorem{pr}[subsection]{Proposition}
\newtheorem{lm}[subsection]{Lemma}

\newtheorem{df}[subsection]{Definition}

\theoremstyle{definition}
\newtheorem{exa}[subsection]{Example}
\theoremstyle{remark}
\newtheorem{rmk}[subsection]{Remark}

\newcommand{\rH}{{\rm H}}

\newcommand{\sO}{\mathscr{O}}

\newcommand{\Z}{\mathbb{Z}}

\newcommand{\G}{\mathbb{G}}

\newcommand{\proj}{\mathbb P}

\newcommand{\Gal}{\operatorname{Gal}}

\newcommand{\Spec}{\operatorname{Spec}}

\newcommand{\GL}{\operatorname{GL}}
\newcommand{\Sym}{\operatorname{Sym}}

\newcommand{\Ind}{\operatorname{Ind}}
\newcommand{\Ker}{\operatorname{Ker}}

\newcommand{\im}{\operatorname{Im}}

\newcommand{\gen}[1]{\langle #1 \rangle}
\newcommand{\id}{\mathds{1}}
\newcommand{\ideal}[1]{\mathfrak{#1}}
\newcommand{\inv}{^{-1}}

\newcommand{\hidden}[1]{\footnote{Hidden:  #1}}
\renewcommand{\hidden}[1]{}

\renewcommand{\phi}{\varphi}

\begin{document}
\pagestyle{plain}
\title{Splitting varieties for cup products with $\Z/3$-coefficients}

\author{Brandon Boggess}\thanks{This work is supported by the REU supplement to NSF DMS-1406380 of PI Kirsten Wickelgren}
\address{School of Mathematics, Georgia Institute of Technology, Atlanta~GA}
\email{bboggess3@gatech.edu}
\date{\today}
\begin{abstract}
We connect Veronese embeddings to splitting varieties of cup products. We then
give an algorithm for constructing splitting varieties for cup products with
$\Z/n$ coefficients, with an explicit calculation for $n=3$. An application to
the automatic realization of Galois groups is given.
\end{abstract}
\maketitle

{\parskip=12pt 

\section{Introduction}

For a functorial assignment $\eta$ of a cohomology class $\eta_F \in \rH^*(\Spec F, \Z/n)$ to fields $F$ over some ground field $k$,
a splitting variety is a scheme $X$ over $k$ which has $F$-points if and only if $\eta_F$ vanishes.
Let $k$ be a number field containing an $n$th root of unity and $F$ a field extension of $k$.
In this paper, we give an algorithm for constructing splitting varieties for cup
products of elements of $\rH^1(\Spec F,\Z/n)$, with an explicit calculation for
$n=3$.

It has been known for a while that quotient schemes give rise to splitting varieties. In
particular, fixed points give rise to versal torsors, and we can use versal torsors
to construct splitting varieties \cite[Example 5.4]{GMS} \cite[Prop. 4.11]{BF}. 
We will construct an algorithm to compute such
fixed points and then run it. What comes out of this run through is long and complicated,
but this is not surprising -- the key ingredient in the algorithm is the Veronese
embedding, and this is known to be complicated for high degree and dimension.
To explain it, we must first develop some notation.

Let $H$ be the group of upper triangular $3 \times 3$-matrices with diagonal entries all 1 and coefficients in $\Z/n$; this is the mod n Heisenberg group.
Let $a_{ij}: H \to \Z/(n)$ be the function taking a matrix to its $(i,j)$-entry. Denote by $E_{ij}$ the matrix such that $a_{ij}(E_{ij}) = 1$ and $a_{kl}(E_{ij}) = 0$ for $k \ne l$.

Let $N \subset H$ be the subgroup $N = \Ker(a_{12}: H \to \Z/n)$.
We have a 1-dimensional representation of $N$ via the map $\rho: N \to \GL_1(k)$ defined by $g \mapsto \zeta_n^{a_{13}(g)}$, where $\zeta_n = e^{2 \pi i/n}$.

We first define a representation $V \cong k^n$ of $H$ by the induced representation
$\Ind_N^H \rho$. From this, we define another representation
\begin{equation*}
  \sigma = \Ind_N^H \rho \times a_{12} \times a_{23} : H \to \GL(V) \times \GL_2(k) ,
\end{equation*}
of $H$, where $a_{12}$ and $a_{23}$ are one dimensional representations given by $g \mapsto \zeta_n^{a_{12}(g)}$ and $g \mapsto \zeta_n^{a_{23}(g)}$, respectively. 
Finding the fixed ring of $k[V] \otimes_k k[\alpha,\alpha\inv ,\beta,\beta\inv ]$ under the representation $\sigma$
gives our splitting variety.

Since $E_{13}$ acts by $x_i \mapsto \zeta_n x_i$ and fixes $\alpha^{\pm 1}$ and $\beta^{\pm1}$, we find that the fixed ring
\[
  k[x_1,\dotsc,x_n,\alpha,\alpha\inv ,\beta,\beta\inv ]^{\gen{E_{13}}} \cong k[\Sym^n V] \otimes_k k[\alpha,\alpha\inv ,\beta,\beta\inv ].
\]
Since the quotient $H/\gen{E_{13}} \cong \Z/n \times \Z/n$ acts on $k[\Sym^n V]$, there is a 
basis $\{E_{12},E_{23}\}$ of $\Z/n \times \Z/n$ and a $\Z/n \times \Z/n$- action on $\Sym^n V$.
But $E_{12}$ and $E_{23}$ commute in $H/\gen{E_{13}}$, so we can find a simultaneous eigenbasis of $\Sym^n V$ with respect to the actions of $E_{12}$ and $E_{23}$.
Let $N = \binom{2n-1}{n}$ denote the degree of $\Sym^n V$. We will define a surjection
\[
  \theta: k[z_1,\dotsc,z_N,\alpha,\alpha\inv ,\beta,\beta\inv ] \to k[x_1,\dotsc,x_n,\alpha,\alpha\inv ,\beta,\beta\inv ]^H
\]
by taking the $z_i$ to these eigenvectors weighted into the $(1,1)$-eigenspace -- since $\alpha$ and $\beta$ have eigenvalues $(\zeta_n,1)$ and $(1,\zeta_n)$ under the
action of $(E_{12},E_{23})$, respectively, this can be achieved by multiplying by powers of $\alpha$ and $\beta$.

It remains to find the kernel of $\theta$.
To do this, we define a map $\pi: k[w_1,\dotsc,w_N] \to k[\Sym^n V]$ by taking generators to distinct degree $n$ monomials. The kernel of this map cuts out the toric ideal of the $n$th Veronese embedding
of $\proj^{n-1}$ \cite[14.1]{Stu}. This ideal has been widely studied, and Sturmfels gives an algorithm to compute a Gr\"obner basis for it. This gives a basis of
$\Sym^n V$, and we can define a change of basis map $p$ from this basis to the eigenbasis under the action of $\Z/n \times \Z/n$.

We prove (Proposition \ref{pr:kernel_theta}) that one can take a known generating set for the toric ideal of the $n$th Veronese embedding of $\proj^{n-1}$ and
obtain the kernel of $\theta$. The proof of this does not depend on $n$ being 3, thus giving a way of algorithmically computing splitting varities of
cup products with coefficients in $\Z/n$.

Finally, we will prove that this construction indeed produces a splitting variety.
In order to state this result, we need a few more definitions.
Let $\kappa: F^* \to \rH^1(\Spec F, \Z/3)$ be the Kummer map, i.e. the map obtained by applying $\rH^*(\Spec F, -)$ to the short exact sequence
\[
  1 \to \mu_3 \to \G_m \xrightarrow{z \mapsto z^3} \G_m \to 1
\]
and identifying $\mu_3 \cong \Z/3$. Let $a,b \in F^*$. We perform a lengthy computation
leading to the main theorem, and so refer to the result of that computation
rather than restating it: let $X(a,b)$ be defined as in Definition
\ref{df:splitting_variety}.
\begin{tm}\label{tm:main_thm_intro}
The scheme $X(a,b)$ has an $F$-point if and only if $\kappa(a) \smile \kappa(b) = 0$ in $\rH^2(\Spec F,\Z/3)$.
\end{tm}

Theorem \ref{tm:main_thm_intro} gives the following automatic realization result for Galois
groups, which we will discuss in Section \ref{sec:galois}.
\begin{tm}
  Suppose $a,b\in F^*$ are such that $F(\sqrt[3] a,\sqrt[3] b)/F$ is a
  $\Z/3\times\Z/3$-Galois extension. Then the following are equivalent:
  \begin{enumerate}
    \item There exists a $\Z/3$-Galois extension $L/F(\sqrt[3] a,\sqrt[3] b)$
      such that $L/F$ is an $H$-Galois extension;
    \item $\kappa(a)\smile\kappa(b)=0$ in $\rH^2(\Spec F,\Z/3)$;
    \item The scheme $X(a,b)$ has an $F$-point.
  \end{enumerate}
\end{tm}

\renewcommand{\abstractname}{Acknowledgements}
\begin{abstract}
I would like to thank Kirsten Wickelgren for introducing me to this problem and sharing ideas. I also thank Josephine Yu for showing me \cite{Stu} and our discussion about Veronese embeddings.
\end{abstract}

\section{Splitting Variety}\label{sec:splitting_variety}

Let $H$ be the group of upper triangular $3 \times 3$-matrices with diagonal entries all 1 and coefficients in $\Z/3$; this is the mod 3 Heisenberg group.
Let $a_{ij}: H \to \Z/(3)$ be the function taking a matrix to its $(i,j)$-entry. Denote by $E_{ij}$ the matrix such that $a_{ij}(E_{ij}) = 1$ and $a_{kl}(E_{ij}) = 0$ for $k \ne l$.

Let $N \subset H$ be the subgroup $N = \Ker(a_{12}: H \to \Z/3)$. We have a 1-dimensional representation of $N$ via the map $\rho: N \to \GL_1(k)$ defined by $g \mapsto \zeta^{a_{13}(g)}$, where $\zeta = e^{2 \pi i/3}$.

We get an induced representation $\Ind_N^H \rho$ of $H$. Following the procedure in \cite[pp.~32-33]{FH}, we can find a homomorphism $H \to \GL(V) \cong \GL_3(k)$ corresponding to $\Ind_N^H \rho$. Given a basis
$\{u_0,u_1,u_2\}$ of $V$ and $g \in H$, this representation acts on $V$ by
\[
  g u_i = \zeta^{a_{13}(g) - (i + a_{12}(g))a_{13}(g)} u_{i + a_{12}(g) \bmod 3}.
\]

Let
\begin{equation} \label{eq:representation_sigma}
  \sigma = \Ind_N^H \rho \times a_{12} \times a_{23} : H \to \GL(V) \times \GL_2(k) ,
\end{equation}
where $a_{12}$ and $a_{23}$ are one dimensional representations given by $g \mapsto \zeta^{a_{12}(g)}$ and $g \mapsto \zeta^{a_{23}(g)}$, respectively. 
This gives another representation of $H$.

Let
\[
  k[x_1,\dotsc,x_n]^H = \{f \in k[x_1,\dotsc,x_n] : gf =f  \quad \forall g \in H \}
\]

We would like to express $k[x_1,x_2,x_3,\alpha,\alpha\inv ,\beta,\beta\inv ]^H$ as a quotient $k[z_1,\dotsc,z_N]/\gen{f_1,\dotsc,f_m}$ for the representation
$\sigma$ of \eqref{eq:representation_sigma}. Define a map
\[
  \theta: k[z_1,\dotsc,z_{10},a,a\inv ,b,b\inv ] \to k[x_1,x_2,x_3,\alpha,\alpha\inv ,\beta,\beta\inv ]
\]
by
\begin{align} \label{eq:theta_def}
  z_1 &\mapsto x_1^3 + x_2^3 + x_3^3 & z_6 &\mapsto \alpha \beta^2 (\zeta x_1^2x_3 + \zeta^2 x_2^2x_1 + x_3^2x_2) & a &\mapsto \alpha^3 \\
  z_2 &\mapsto \alpha^2 (\zeta^2 x_1^3 + \zeta x_2^3 + x_3^3) & z_7 &\mapsto \beta (x_1^2x_2 + x_2^2x_3 + x_3^2x_1) & a\inv  &\mapsto \alpha^{-3} \nonumber \\
  z_3 &\mapsto \alpha (\zeta x_1^3 + \zeta^2 x_2^3 + x_3^3) & z_8 &\mapsto \alpha^2 \beta (\zeta^2 x_1^2x_2 + \zeta x_2^2x_3 + x_3^2x_1) & b &\mapsto \beta^3\nonumber \\
  z_4 &\mapsto \beta^2 (x_1^2x_3 + x_2^2x_1 + x_3^2x_2) & z_9 &\mapsto \alpha \beta (\zeta x_1^2x_2 + \zeta^2 x_2^2x_3 + x_3^2x_1) & b\inv  &\mapsto \beta^{-3} \nonumber \\
  z_5 &\mapsto \alpha^2 \beta^2 (\zeta^2 x_1^2x_3 + \zeta x_2^2x_1 + x_3^2x_2) & z_{10} &\mapsto x_1x_2x_3 & \nonumber
\end{align}
The following proposition gives an explanation for this map.

\begin{pr}\label{pr:theta_surjection}
The image of $\theta$ is $k[x_1,x_2,x_3, \alpha, \alpha\inv , \beta, \beta\inv ]^H$.
\end{pr}

\bpf
  The action of $E_{13}$ takes $E_{13}(x_i) = \zeta x_i$ for $i=1,2,3$ and fixes $\alpha^{\pm1}$ and $\beta^{\pm1}$. Thus
  \[
    k[x_1,x_2,x_3,\alpha,\alpha\inv ,\beta,\beta\inv ]^{\gen{E_{13}}} \cong k[\Sym^3 V] \otimes_k k[\alpha,\alpha\inv ,\beta,\beta\inv ]
  \]
  Now, $H/\gen{E_{13}} \cong \Z/3 \times \Z/3$ has a $\Z/3$ basis $\{E_{12},E_{23}\}$.
  The action of $E_{12}$ permutes the $x_i$ by $(x_1\;x_2\;x_3)$ (here we use cycle notation, cf. \cite[$\S$ 1.5]{Art}) and takes $E_{12}(\alpha^{\pm1}) = \zeta \alpha^{\pm1}, E_{12}(\beta^{\pm1}) = \beta^{\pm1}$.
  The action of $E_{23}$ is given by $E_{23}(x_1) = x_1, E_{23}(x_2) = \zeta^2x_2, E_{23}(x_3) = \zeta x_3, E_{23}(\alpha^{\pm1}) = \alpha^{\pm1}, E_{23}(\beta^{\pm1}) = \beta^{\pm1}$.
  Hence there is a simultaneous eigenbasis for $\Z/3 \times \Z/3$ of $\Sym^3 V \oplus k\alpha \oplus k\alpha\inv  \oplus k\beta \oplus k\beta\inv $ given by
  \begin{equation} \label{eq:eigenvectors}
    \begin{array}{c c}
      \text{Eigenvector} & \text{Eigenvalues} \\
      x_1^2 x_3 + x_2^2 x_1 + x_3^2 x_2 & (1, \zeta)  \\
      x_1^2 x_2 + x_2^2 x_3 + x_3^2 x_1 & (1, \zeta^2)  \\
      x_1^3 + x_2^3 + x_3^3             & (1,1)   \\
      x_1 x_2 x_3                       & (1,1)   \\
      \zeta^2 x_1^2 x_3 + \zeta x_2^2 x_1 + x_3^2 x_2 & (\zeta, \zeta)  \\
      \zeta^2 x_1^2 x_2 + \zeta x_2^2 x_3 + x_3^2 x_1 & (\zeta, \zeta^2)  \\
      \zeta^2 x_1^3 + \zeta x_2^3 + x_3^3             & (\zeta, 1)    \\
      \zeta x_1^2 x_3 + \zeta^2 x_2^2 x_1 + x_3^2 x_2 & (\zeta^2, \zeta)  \\
      \zeta x_1^2 x_2 + \zeta^2 x_2^2 x_3 + x_3^2 x_1 & (\zeta^2, \zeta^2)  \\
      \zeta x_1^3 + \zeta^2 x_2^3 + x_3^3             & (\zeta^2, 1)  \\
      \alpha & (\zeta,1) \\
      \alpha\inv  & (\zeta^2,1) \\
      \beta  & (1,\zeta) \\
      \beta\inv   & (1,\zeta^2)
    \end{array}
  \end{equation}

  Label these eigenvectors as $v_1,\dotsc,v_{14}$. Then any
  element $f \in k[\Sym^3 V] \otimes_k k[\alpha,\alpha\inv ,\beta,\beta\inv ]$ can be expressed as a polynomial in the $v_i$. Since each monomial
  of $f$ in the variables $v_i, \alpha, \beta$ goes to a scalar multiple of itself under the action of any element in $H$,
  $f$ is fixed by $H$ if and only if each monomial is. Therefore $k[x_1,x_2,x_3,\alpha,\alpha\inv ,\beta,\beta\inv ]^H$
  is generated by products of the eigenvectors whose eigenvalues under $E_{12}$ and $E_{23}$
  multiply to 1, which we claim is the image of $\theta$.

  Let
  $(\zeta^{m_i},\zeta^{n_i})$ be the eigenvalues of $v_i$ under $(E_{12},E_{23})$, respectively. Then, re-indexing the $v_i$ as necessary,
  \[
    \theta(z_i) = \alpha^{3-m_i} \beta^{3-n_i} v_i
  \]
  by \eqref{eq:theta_def} (this even works for $\alpha^{\pm1}$ and $\beta^{\pm1}$,
  considering $z_{11}$ through $z_{14}$ to be $a,a\inv ,b,b\inv $, respectively).
  If the monomial $v_{1}^{r_1}\cdots v_{14}^{r_14}$ is fixed, then $\sum r_i m_i$ and
  $\sum r_i n_i$ are both divisible by 3, so that the eigenvalue of the product is $(1,1)$. Hence
  \[
    \prod_{i=1}^{14} \theta(z_{i}^{r_i}) = \prod_{i=1}^{14} \alpha^{r_i(3-m_i)} \beta^{r_i(3-n_i)} v_{i}^{r_i}
      = \alpha^{3m} \beta^{3n} \prod_{i=1}^{14} v_{i}^{r_i} .
  \]
  for some $m,n \in \Z$. Therefore
  \[
    \prod_{i=1}^{14} v_{i}^{r_i}
      = \theta(a^{-m}) \theta(b^{-n}) \prod_{i=1}^{14} \theta(z_{i}^{r_i}) \in \im(\theta) .
  \]
\epf

The eigenspace decomposition of $k[x_1,x_2,x_3,\alpha,\alpha\inv ,\beta,\beta\inv ]^{\gen{E_{13}}}$ means that any degree 3 monomial
in $k[x_1,x_2,x_3]$ can be expressed in terms of the eigenvectors for $\Z/3 \times \Z/3$. This will be useful in finding the kernel of $\theta$.

Define the map $\pi: k[w_1,\dots,w_{10}] \to k[\Sym^3 V]$ by $w_1 \mapsto x_1x_2x_3$, $w_2 \mapsto x_1^2$, etc.
The kernel of $\pi$ cuts out the toric ideal of the 3rd Veronese embedding of $\proj^2$ \cite{Stu}, and \texttt{Macaulay2} \cite{M2} can find the generators:
\begin{align} \label{eq:kernel_pi}
  \Ker \pi = \langle &w_8w_9-w_1w_{10}, w_7w_9-w_5w_{10}, w_5w_9-w_2w_{10}, w_3w_9-w_7w_{10}, w_1w_9-w_6w_{10}, \\
  &\quad w_8^2-w_3w_{10}, w_6w_8-w_5w_{10}, w_4w_8-w_{10}^2, w_1w_8-w_7w_{10}, w_6w_7 - w_2w_8, \nonumber \\
  &\quad w_4w_7-w_1w_{10}, w_1w_7-w_5w_8, w_6^2-w_2w_9, w_4w_6-w_9^2, w_3w_6-w_5w_8, \nonumber \\
  &\quad w_1w_6-w_2w_{10}, w_5^2-w_2w_7, w_4w_5-w_6w_{10}, w_3w_5-w_7^2, w_1w_5-w_2w_8, \nonumber \\
  &\quad w_3w_4-w_8w_{10}, w_2w_4-w_6w_9, w_1w_4-w_9w_{10}, w_2w_3-w_5w_7, w_1w_3-w_7w_8,\nonumber \\
  &\quad w_1w_2-w_5w_6, w_1^2-w_5w_{10} \nonumber \rangle.
\end{align}
There are 27 generators, each of which falls into one of five categories:
\begin{enumerate}[1)]
  \item Relations of the form $x_i^2 x_j \cdot x_j^2 x_i = x_i^3 \cdot x_j^3$. There are three such relations.
    These are mapped to by $w_5w_7 - w_2w_3$, $w_6w_9 - w_2w_4$, and $w_8w_{10}-w_3w_4$.

  \item Relations of the form $x_i^2 x_j \cdot x_i^2 x_k = x_i^3 \cdot x_ix_jx_k$. There are three of these relations.
    These are mapped to by $w_5w_6 - w_1w_2$, $w_7w_8 - w_1w_3$, and $w_9w_{10} - w_1w_4$.

  \item Relations of the form $(x_i^2 x_j)^2 = x_i^3 \cdot x_j^2 x_i$. There are six of these, and they are mapped to by
    $w_5^2 - w_2w_7$, $w_6^2 - w_2w_9$, $w_7^2 - w_3 w_5$, $x_8^2 - w_3w_{10}$, $w_9^2 - w_4w_6$, and $w_{10}^2 - w_4w_8$.

  \item Relations of the form $x_i^2x_j \cdot x_j^2 x_k = x_j^3 \cdot x_i^2 x_k = x_j^2 x_i \cdot x_ix_jx_k$. There are six products of this
    form, each of which gives two relations; thus there are 12 relations of this form.

  \item There are three ways to write $(x_ix_jx_k)^2$ in terms of the other monomials.
    These are mapped to by $w_1^2 - w_5w_{10}$, $w_1^2 - w_6w_8$, and $w_1^2 - w_7w_9$.
\end{enumerate}
These give exactly all 27 generators for $\Ker \pi$.

\begin{rmk} \label{rmk:extension_pi}
Extending $\pi$ to the map
\[
  \pi \otimes \id: k[w_1,\dots,w_{10}] \otimes_k k[\alpha,\alpha\inv ,\beta,\beta\inv ] \to k[\Sym^3 V] \otimes_k k[\alpha,\alpha\inv ,\beta,\beta\inv ]
\]
does not affect the kernel because $k[\alpha,\alpha\inv ,\beta,\beta\inv ]$ is a flat $k$-algebra.
\end{rmk}

As in the proof of Proposition \ref{pr:theta_surjection}, we label the eigenvectors of $\Sym^3 V$ as $v_1,\dotsc,v_{10}$ in such a way that $\theta(z_i)$ is a multiple of $v_i$. Let $\pi': k[v_1,\dotsc,v_{10}] \to k[\Sym^3 V]$ take $v_i$ to its expression in \eqref{eq:eigenvectors}.

Since $H$ acts on $V$, the quotient $H/\gen{E_{13}}$ acts on the fixed elements $k[V]^{\gen{E_{13}}} \cong k[\Sym^3 V]$, so there is a $\Z/3 \times \Z/3$ action on $k[\Sym^3 V]$ with basis $\{E_{12},E_{23}\}$.
This gives a natural action on $k[v_1,\dotsc,v_{10}]$. For any $g \in \Z/3 \times \Z/3$, we have that $g \pi'(v_i) = \zeta^{k_i} \pi'(v_i)$ for some $k_i$. We then define
\[
  g v_i := \zeta^{k_i} v_i .
\]
It follows that $\pi'$ is equivariant.

Denote by $p$ the map $k[w_1,\dotsc,w_{10}] \to k[v_1,\dotsc,v_{10}]$ resulting from the change of basis corresponding to the two bases of $\Sym^3 V$ given
by $\{w_1,\dotsc,w_{10}\}$ and $\{v_1,\dotsc,v_{10}\}$. Note that
\[
  \begin{tikzcd}
    k[w_1,\dots,w_{10}] \arrow{r}{p} \arrow{dr}{\pi} &
      k[v_1,\dotsc,v_{10}] \arrow{d}{\pi'} \\
    & k[\Sym^3 V]
  \end{tikzcd}
\]
commutes and that $p$ is an isomorphism. We want a map $\phi: k[z_1,\cdots,z_{10},a,a\inv ,b,b\inv ] \to k[w_1,\dotsc,w_{10},\alpha,\alpha\inv ,\beta,\beta\inv ]$ 
 such that tensoring with $k[\alpha,\alpha\inv ,\beta,\beta\inv ]$ makes
\begin{equation} \label{eq:diagram_all_maps}
  \begin{tikzcd}
    k[z_1,\dotsc,z_{10},\alpha,\alpha\inv ,\beta,\beta\inv ] \arrow{r}{\phi} \arrow{dr}{\theta} & k[w_1,\dotsc,w_{10},\alpha,\alpha\inv ,\beta,\beta\inv ] \arrow{d}{\pi \otimes \id} \arrow{r}{p \otimes \id} & k[v_1,\dotsc,v_{10},\alpha,\alpha\inv ,\beta,\beta\inv ] \arrow{dl}{\pi' \otimes \id} \\
    & k[\Sym^3 V] \otimes_k k[\alpha,\alpha\inv ,\beta,\beta\inv ] \arrow[hook]{d}& \\
    & k[V] \otimes_k k[\alpha,\alpha\inv ,\beta,\beta\inv ] &
  \end{tikzcd}
\end{equation}
commute. This map is given by
\begin{align} \label{eq:phi_def}
    z_1 &\mapsto w_2 + w_3 + w_4 &z_6 &\mapsto \alpha \beta^2(\zeta w_6 + \zeta^2 w_7 + w_{10}) &a &\mapsto \alpha^3 \\
    z_2 &\mapsto \alpha^2 (\zeta^2 w_2 + \zeta w_3 + w_4) &z_7 &\mapsto \beta(w_5 + w_8 + w_9) &a\inv  &\mapsto \alpha^{-3} \nonumber \\
    z_3 &\mapsto \alpha(\zeta w_2 + \zeta^2 w_3 + w_4) &z_8 &\mapsto \alpha^2 \beta(\zeta^2 w_5 + \zeta w_8 + w_9) &b &\mapsto \beta^{3} \nonumber \\
    z_4 &\mapsto \beta^2(w_6 + w_7 + w_{10}) &z_9 &\mapsto \alpha\beta(\zeta w_5 + \zeta^2 w_8 + w_9) &b\inv  &\mapsto \beta^{-3} \nonumber \\
    z_5 &\mapsto \alpha^2\beta^2(\zeta^2 w_6 + \zeta w_7 + w_{10}) &z_{10} &\mapsto w_1 \nonumber
\end{align}
as can be verified by comparing \eqref{eq:phi_def} to \eqref{eq:theta_def}.
The map $\phi$ is used to bring results about Veronese embeddings into the picture, and $p$ allows us to translate this information into statements about eigenspaces under $\Z/3 \times \Z/3$.

Given a generator $w_iw_j - w_kw_l$ of $\Ker \pi$, we can construct an element of $\Ker \theta$ as follows: let $h = p(w_iw_j - w_kw_l)$.
Obeserve that $\pi'(h) = 0$. Each monomial in $h$ belongs to some eigenspace; let $h_{m,n}$ be the projection of $h$ onto the $(\zeta^m,\zeta^n)$ eigenspace, so that
\[
  h = \sum_{0 \le m,n < 3} h_{m,n} .
\]
Because $\pi'$ is equivariant, the $\pi'(h_{m,n})$ are in distinct eigenspaces. It follows that $\pi'(h)=0$ if and only if each $\pi'(h_{m,n}) = 0$.
By construction, each $\alpha^{-m} \beta^{-n} \pi'(h_{m,n}) \in \im \theta$, with the exponents of $\alpha$ and $\beta$ taken mod 3. Thus if $\pi'(h) = 0$,
the elements $\alpha^{-m} \beta^{-n} \pi'(h_{m,n})$ are in $\Ker \theta$.

\begin{exa}
We will construct elements of $\Ker \theta$ from $w_8w_9 - w_1w_{10} \in \Ker \pi$ using this method. We first apply $p$ to get
\[
  h = \frac{1}{9} (v_7 + \zeta^2 v_8 + \zeta v_9)(v_7 + v_8 + v_9) - \frac{1}{3}v_{10}(v_4 + v_5 + v_6) \in k[v_1,\dotsc,v_{10]}]
\]
We get the decomposition into eigenspaces (ignoring those with no contribution)
\begin{align*}
  h_{0,1} &= v_7^2 - v_8v_9 - 3v_4v_{10} \\
  h_{1,1} &= \zeta v_9^2 - \zeta v_7v_8 - 3v_5v_{10} \\
  h_{2,1} &= \zeta^2 v_8^2 - \zeta^2 v_7v_9 - 3v_6v_{10} .
\end{align*}
We then see from \eqref{eq:theta_def} that
\begin{align*}
  \theta\left(z_7^2 - \frac{1}{a}z_8z_9 - 3 z_4z_{10}\right) &= \beta^2 \pi'(h_{0,1}) \\
  \theta\left(\zeta z_9^2 - \zeta z_7z_8 - 3 z_5z_{10}\right) &= \alpha^2 \beta^2 \pi'(h_{1,1}) \\
  \theta\left(\frac{\zeta^2}{a} z_8^2 - \zeta^2 z_7z_9 - 3 z_6z_{10}\right) &= \alpha \beta^2 \pi'(h_{2,1}) ,
\end{align*}
and it is straightforward to verify that these elements are indeed in $\Ker \theta$.
\end{exa}

We now show that all elements of $\Ker \theta$ come about from this process.

\begin{pr} \label{pr:kernel_theta}
Given $\gen{u_1,\dotsc,u_M}$ a generating set for $\Ker \pi$, we can find a set of generators for $\Ker \theta$.
\end{pr}

\bpf
Since the maps $\phi$ and $p \otimes \id$ are inejctive, we can view $k[z_1,\dotsc,z_{10},\alpha,\alpha\inv ,\beta,\beta\inv ]$ as a subring of
$k[v_1,\dotsc,v_{10},\alpha,\alpha\inv ,\beta,\beta\inv ]$. By Proposition \ref{pr:theta_surjection}, this subring is the $(1,1)$-eigenspace. In view of
\eqref{eq:diagram_all_maps}, finding $\Ker(\theta)$ can be reduced to finding the intersection of $\Ker(\pi' \otimes \id)$ with the $(1,1)$-eigenspace of
$k[v_1,\dotsc,v_{10},\alpha,\alpha\inv ,\beta,\beta\inv ]$.

Given generators for $\Ker(\pi)$, we can obtain generators for $\Ker(\pi' \otimes \id)$ as follows. By Remark \ref{rmk:extension_pi}, we then know that $\Ker(\pi \otimes \id)= \Ker(\pi) \otimes_k k[\alpha,\alpha\inv ,\beta,\beta\inv ]$.
Then, since $p \otimes \id$ is an isomorphism, we find that $\Ker(\pi' \otimes \id)= p \otimes \id(\Ker(\pi \otimes \id))$. In particular, given finitely
many generators $u_1,\dotsc,u_M$ of $\Ker \pi$, we get that $\Ker(\pi' \otimes \id) = \gen{p(u_1) \otimes 1, \dotsc, p(u_M) \otimes 1}$.

Let $h_1, \dotsc, h_N$ denote the projections onto eigenspaces of the generators of $\ideal a = \Ker(\pi' \otimes \id)$. Note that the ring
$k[v_1,\dotsc,v_{10},\alpha,\alpha\inv ,\beta,\beta\inv ]$ is a direct sum of eigenspaces under the action of $\Z/3 \times \Z/3$. That $\ideal{a} \subset \gen{h_1,\dotsc,h_N}$
is immediate, and the other direction follows because $\pi'$ is equivariant. Indeed, if $h \in \ideal a$ and $h_{m,n}$ is its projection onto the
$(\zeta^m,\zeta^n)$-eigenspace, then
\[
  \sum_{0 \le m,n < 3} \pi' \otimes 1 (h_{m,n}) = 0.
\]
But each of these terms is in a distinct eigenspace because $\pi'$ is equivariant,
and so each term must itself be 0.

Let $h_i$ be in the $(\zeta^{m_i},\zeta^{n_i})$-eigenspace, so that $\alpha^{-m_i}\beta^{-n_i}h_i$ is in the $(1,1)$-eigenspace. We claim that $\gen{\alpha^{-m_1}\beta^{-n_1}h_1,\dotsc,\alpha^{-m_N}\beta^{-n_N}h_N}$ is the intersection of $\ideal a$ with the $(1,1)$-eigenspace,
which we will denote by $\ideal{a}_{(1,1)}$. The inclusion $\gen{\alpha^{-m_1}\beta^{-n_1}h_1,\dotsc,\alpha^{-m_N}\beta^{-n_N}h_N} \subset \ideal{a}_{(1,1)}$ follows immediately.
Now let $f \in \ideal{a}_{(1,1)}$. Then we can write
\[
  f = \sum_{i=1}^N g_i h_i
\]
where each $g_i \in k[v_1,\dotsc,v_{10},\alpha,\alpha\inv ,\beta,\beta\inv ]$. In order for $f$ to be in the $(1,1)$-eigenspace, all of the terms $g_i h_i$ in
a different eigenspace must cancel out; hence we may assume that each $g_i h_i$ is in the $(1,1)$-eigenspace. By assumption, $h_i$ is in the $(\zeta^{m_i},\zeta^{n_i})$-eigenspace,
so $g_i$ is in the $(\zeta^{-m_i},\zeta^{-n_i})$-eigenspace. Let $g_i' = \alpha^{m_i} \beta^{n_i} g_i$. Then $g_i'$ is in the $(1,1)$-eigenspace, and
\[
  f = \sum_i g_i' (\alpha^{-m_i}\beta^{-n_i}h_i).
\]
\epf

\begin{rmk}
Given generators for the toric ideal of the appropiate Veronese embedding, Proposition \ref{pr:kernel_theta} gives a way to algorithmically compute splitting varieties for cup products with coefficients in $\Z/n$ for any $n$.
\end{rmk}

Performing this calculation for each generator in \eqref{eq:kernel_pi} yields
\begin{pr}\label{pr:theta_kernel}
The kernel of $\theta$ is generated by
\begin{itemize}
  \item $z_7^2 - \frac{1}{a} z_8 z_9 - 3 z_4 z_{10}$
  \item $\frac{\zeta^2}{a} z_8^2 - \zeta^2 z_7 z_9 - 3 z_6 z_{10}$
  \item $\zeta z_9^2 - \zeta z_7 z_8 - 3 z_5z_{10}$

  \item $(1-\zeta)z_4z_8 + (\zeta^2-1)z_5z_7 + (\zeta-\zeta^2)z_6z_9$
  \item $(1-\zeta^2)z_4z_9 + \frac{\zeta^2-\zeta}{a}z_5z_8 + (\zeta-1)z_6z_7$

  \item $z_7^2 - \frac{1}{a}z_8z_9 - z_1z_4 - \frac{\zeta}{a}z_2z_6 - \frac{\zeta^2}{a}z_3z_5$
  \item $\frac{\zeta}{a}z_8^2 - \zeta z_7z_9 - z_1z_6 - \frac{\zeta}{a}z_2z_5 - \zeta^2 z_3z_4$
  \item $\zeta^2 z_9^2 -\zeta^2 z_7z_8 - z_1z_5 - \zeta z_2z_4 - \zeta^2 z_3z_6$

  \item $z_1z_7 + \frac{\zeta^2}{a}z_2z_9 + \frac{\zeta}{a}z_3z_8 - \frac{1}{b}z_4^2 + \frac{1}{ab}z_5z_6$
  \item $z_1z_8 + \zeta^2 z_2z_7 + \zeta z_3z_9 - \frac{\zeta}{b}z_6^2 + \frac{\zeta}{b}z_4z_5$
  \item $z_1z_9 + \frac{\zeta^2}{a}z_2z_8 + \zeta z_3z_7 - \frac{\zeta^2}{ab}z_5^2 + \frac{\zeta^2}{b}z_4z_6$

  \item $\frac{1}{b}z_4^2 - \frac{1}{ab}z_5z_6 - 3 z_7z_{10}$
  \item $\frac{\zeta}{ab}z_5^2 - \frac{\zeta}{b}z_4z_6 - 3z_9z_{10}$
  \item $\frac{\zeta^2}{b}z_6^2 - \frac{\zeta^2}{b}z_4z_5- 3z_8z_{10}$

  \item $z_7^2 + \frac{2}{a}z_8z_9 - z_1z_4 - \frac{\zeta^2}{a}z_2z_6 -\frac{\zeta}{a}z_3z_5$
  \item $\frac{\zeta}{a}z_8^2 + 2\zeta z_7z_9 - z_1z_6 -\frac{\zeta^2}{a}z_2z_5 - \zeta z_3z_4 $
  \item $\zeta^2 z_9^2 + 2\zeta^2 z_7z_8 - z_1z_5 - \zeta^2 z_2z_4 - \zeta z_3z_6$

  \item $(\zeta^2-\zeta)z_4z_8 + (\zeta-1)z_5z_7 + (1-\zeta^2)z_6z_9$
  \item $(\zeta-\zeta^2)z_4z_9 + \frac{1-\zeta}{a}z_5z_8 + (\zeta^2-1)z_6z_7$

  \item $z_1z_7 + \frac{\zeta}{a}z_2z_9 + \frac{\zeta^2}{a}z_3z_8 - \frac{1}{b} z_4^2 - \frac{2}{ab} z_5z_6$
  \item $\zeta^2 z_1z_8 + z_2z_7 + \zeta z_3z_9 - \frac{1}{b} z_6^2 - \frac{2}{b}z_4z_5$
  \item $\zeta z_1z_9 + \frac{\zeta^2}{a}z_2z_8 + z_3z_7 - \frac{1}{ab} z_5^2 - \frac{2}{b}z_4z_6$

  \item $\frac{1}{b}z_4^2 - \frac{1}{ab}z_5z_6 - 3z_7z_{10}$
  \item $\frac{\zeta^2}{ab}z_5^2 - \frac{\zeta^2}{b}z_4z_6 - 3 \zeta z_9z_{10}$
  \item $\frac{\zeta}{b}z_6^2 - \frac{\zeta}{b}z_4z_5 - 3 \zeta^2 z_8z_{10}$

  \item $\frac{1}{b}z_4^2 - \frac{1}{ab}z_5z_6 - z_1z_7 - \frac{\zeta^2}{a}z_2z_9 - \frac{\zeta}{a}z_3z_8$
  \item $\frac{1}{ab}z_5^2 - \frac{1}{b}z_4z_6 - \zeta z_1z_9 - \frac{1}{a}z_2z_8 - \zeta^2z_3z_7$
  \item $\frac{1}{b}z_6^2 - \frac{1}{b}z_4z_5 - \zeta^2 z_1z_8 - \zeta z_2z_7 - z_3z_9$

  \item $z_1z_4 + \frac{\zeta}{a}z_2z_6 + \frac{\zeta^2}{a}z_3z_5 - 3z_4z_{10}$
  \item $\zeta^2 z_1z_5 + z_2z_4 + \zeta z_3z_6 - 3z_5z_{10}$
  \item $\zeta z_1z_6 + \frac{\zeta^2}{a}z_2z_5 + z_3z_4 - 3z_6z_{10}$

  \item $z_7^2 - \frac{1}{a}z_8z_9 - 3z_4z_{10}$
  \item $\frac{1}{a}z_8^2 - z_7z_9 - 3\zeta z_6z_{10}$
  \item $z_9^2 - z_7z_8 - 3\zeta^2 z_5z_{10}$

  \item $\frac{1}{b}z_4^2 + \frac{2}{ab}z_5z_6 - z_1z_7 - \frac{\zeta}{a}z_2z_9 - \frac{\zeta^2}{a}z_3z_8$
  \item $\frac{\zeta^2}{ab}z_5^2 + \frac{2\zeta^2}{b}z_4z_6 - z_1z_9 - \frac{\zeta}{a}z_2z_8 - \zeta^2 z_3z_7$
  \item $\frac{\zeta}{b}z_6^2 + \frac{2\zeta}{b}z_4z_5 - z_1z_8 - \zeta z_2z_7 - \zeta^2 z_3z_9$

  \item $z_1z_4 + \frac{\zeta^2}{a}z_2z_6 + \frac{\zeta}{a}z_3z_5 - z_7^2 - \frac{2}{a}z_8z_9$
  \item $\zeta z_1z_5 + z_2z_4 + \zeta^2 z_3z_6 - z_9^2 - 2z_7z_8$
  \item $\zeta^2 z_1z_6 + \frac{\zeta}{a}z_2z_5 + z_3z_4 - \frac{1}{a}z_8^2 - 2z_7z_9$

  \item $z_1z_4 + \frac{\zeta}{a}z_2z_6 + \frac{\zeta^2}{a}z_3z_5 - z_7^2 + \frac{1}{a}z_8z_9$
  \item $\zeta z_1z_5 + \zeta^2 z_2z_4 + z_3z_6 - z_9^2 + z_7z_8$
  \item $\zeta^2 z_1z_6 + \frac{1}{a}z_2z_5 + \zeta z_3z_4 - \frac{1}{a}z_8^2 + z_7z_9$

  \item $z_1z_4 + \frac{\zeta}{a}z_2z_6 + \frac{\zeta^2}{a}z_3z_5 - 3z_4z_{10}$
  \item $z_1z_5 + \zeta z_2z_4 + \zeta^2 z_3z_6 - 3\zeta z_5z_{10}$
  \item $z_1z_6 + \frac{\zeta}{a}z_2z_5 + \zeta^2 z_3z_4 - 3\zeta^2 z_6z_{10}$

  \item $z_7^2 + \frac{2}{a}z_8z_9 - z_1z_4 - \frac{\zeta^2}{a}z_2z_6 - \frac{\zeta}{a}z_3z_5$
  \item $\frac{\zeta^2}{a}z_8^2 + 2 \zeta^2 z_7z_9 - \zeta z_1z_6 - \frac{1}{a}z_2z_5 - \zeta^2 z_3z_4$
  \item $\zeta z_9^2 + 2\zeta z_7z_8 - \zeta^2 z_1z_5 - \zeta z_2z_4 - z_3z_6$

  \item $z_1z_7 + \frac{\zeta^2}{a}z_2z_9 + \frac{\zeta}{a}z_3z_8 - \frac{1}{b}z_4^2 + \frac{1}{ab}z_5z_6$
  \item $\zeta z_1z_8 + z_2z_7 + \zeta^2 z_3z_9 - \frac{\zeta^2}{b}z_6^2 + \frac{\zeta^2}{b}z_4z_5$
  \item $\zeta^2 z_1z_9 + \frac{\zeta}{a}z_2z_8 + z_3z_7 - \frac{\zeta}{ab}z_5^2 + \frac{\zeta}{b}z_4z_6$

  \item $z_1z_7 + \frac{\zeta}{a}z_2z_9 + \frac{\zeta^2}{a}z_3z_8 - \frac{1}{b}z_4^2 - \frac{2}{ab}z_5z_6$
  \item $\zeta z_1z_8 + \zeta^2 z_2z_7 + z_3z_9 - \frac{\zeta^2}{b}z_6^2 - \frac{2\zeta^2}{b} z_4z_5$
  \item $\zeta^2 z_1z_9 + \frac{1}{a}z_2z_8 + \zeta z_3z_7 - \frac{\zeta}{ab}z_5^2 - \frac{2\zeta}{b} z_4z_6$

  \item $z_1z_7 + \frac{\zeta^2}{a}z_2z_9 + \frac{\zeta}{a}z_3z_8 - 3z_7z_{10}$
  \item $\zeta^2 z_1z_8 + \zeta z_2z_7 + z_3z_9 - 3\zeta z_8z_{10}$
  \item $\zeta z_1z_9 + \frac{1}{a}z_2z_8 + \zeta^2 z_3z_7 - 3\zeta^2 z_9z_{10}$

  \item $z_1^2 - \frac{1}{a}z_2z_3 - \frac{1}{b}z_4z_7 - \frac{\zeta}{ab}z_5z_9 - \frac{\zeta^2}{ab}z_6z_8$
  \item $\frac{\zeta^2}{a}z_2^2 - \zeta^2 z_1z_3 - \frac{\zeta}{b}z_4z_9 - \frac{\zeta^2}{ab}z_5z_8 - \frac{1}{b}z_6z_7$
  \item $\zeta z_3^2 - \zeta z_1z_2 - \frac{\zeta^2}{b}z_4z_8 - \frac{1}{b}z_5z_7 - \frac{\zeta}{b}z_6z_9$

  \item $z_1^2 - \frac{1}{a}z_2z_3 - \frac{1}{b}z_4z_7 - \frac{\zeta}{ab}z_5z_9 - \frac{\zeta^2}{ab}z_6z_8$
  \item $\frac{\zeta}{a}z_2^2 - \zeta z_1z_3 - \frac{1}{b}z_4z_9 -\frac{\zeta}{ab}z_5z_8 - \frac{\zeta^2}{b}z_6z_7$
  \item $\zeta^2 z_3^2 - \zeta^2 z_1z_2 - \frac{1}{b}z_4z_8 - \frac{\zeta}{b}z_5z_7 - \frac{\zeta^2}{b}z_6z_9$

  \item $\frac{1}{b}z_4z_7 + \frac{1}{ab}z_5z_9 + \frac{1}{ab}z_6z_8 - 3z_1z_{10}$
  \item $\frac{1}{b}z_4z_8 + \frac{1}{b}z_5z_7 + \frac{1}{b}z_6z_9 - 3z_2z_{10}$
  \item $\frac{1}{b}z_4z_9 + \frac{1}{ab}z_5z_8 + \frac{1}{b}z_6z_7 - 3z_3z_{10}$

  \item $z_1^2 - \frac{1}{a}z_2z_3 - \frac{1}{b}z_4z_7 - \frac{\zeta}{ab}z_5z_9 - \frac{\zeta^2}{ab}z_6z_8$
  \item $\frac{1}{a}z_2^2 - z_1z_3 - \frac{\zeta^2}{b}z_4z_9 - \frac{1}{ab}z_5z_8 - \frac{\zeta}{b}z_6z_7$
  \item $z_3^2 - z_1z_2 - \frac{\zeta}{b}z_4z_8 - \frac{\zeta^2}{b}z_5z_7 - \frac{1}{b}z_6z_9$

  \item $\frac{1}{b}z_4z_7 + \frac{1}{ab}z_5z_9 + \frac{1}{ab}z_6z_8 - 3z_1z_{10}$
  \item $\frac{\zeta^2}{b}z_4z_8 + \frac{\zeta^2}{b}z_5z_7 + \frac{\zeta^2}{b}z_6z_9 - 3\zeta^2 z_2z_{10}$
  \item $\frac{\zeta}{b}z_4z_9 + \frac{\zeta}{ab}z_5z_8 + \frac{\zeta}{b}z_6z_7 - 3\zeta z_3z_{10}$

  \item $\frac{1}{b}z_4z_7 + \frac{1}{ab}z_5z_9 + \frac{1}{ab}z_6z_8 - 3 z_1 z_{10}$
  \item $\frac{\zeta}{b}z_4z_8 + \frac{\zeta}{b}z_5z_7 + \frac{\zeta}{b}z_6z_9 - 3 \zeta z_2 z_{10}$
  \item $\frac{\zeta^2}{b}z_4z_9 + \frac{\zeta^2}{ab}z_5z_8 + \frac{\zeta^2}{b}z_6z_7 - 3 \zeta^2 z_3z_{10}$

  \item $\frac{1}{b}z_4z_7 + \frac{\zeta^2}{ab}z_5z_9 + \frac{\zeta}{ab}z_6z_8 - 9z_{10}^2$
  \item $\frac{\zeta}{b}z_4z_8 + \frac{1}{b}z_5z_7 + \frac{\zeta^2}{b}z_6z_9$
  \item $\frac{\zeta^2}{b}z_4z_9 + \frac{\zeta}{ab}z_5z_8 + \frac{1}{b}z_6z_7$
\end{itemize}
\end{pr}

Let $I$ be the ideal of $k[z_1,\dotsc,z_{10},a,a\inv ,b,b\inv ]$ generated by the elements listed in Proposition \ref{pr:theta_kernel}. Propositions \ref{pr:theta_surjection} and \ref{pr:theta_kernel} combine to give the following result.
\begin{lm} \label{lm:fixed_ring}
There is an isomorphism
\[
  (k[V]\otimes_k k[\alpha,\alpha\inv\beta,\beta\inv])^H \cong k[z_1,\dotsc,z_{10},a,a\inv ,b,b\inv ] / I.
\]
\end{lm}

Fix $a,b \in k^*$, so that $k[a,a\inv ,b,b\inv ] = k$; we are plugging in values
for $a$ and $b$ in the ring of Lemma \ref{lm:fixed_ring}. 
Define $S$ to be closed subscheme corresponding to
the union of the fixed subspaces of the representation
$\sigma$ of $H$.

\begin{df}\label{df:splitting_variety}
  Define an open subset of $\Spec k[z_1,\dotsc,z_{10}] / I$ by
  \[
    X(a,b) := \Spec k[z_1,\dotsc,z_{10}] / I - S.
  \]
\end{df}
\begin{rmk}
  Subtracting $S$ causes $H$ to act freely on $X(a,b)$.
\end{rmk}

By construction, there is an $H$-torsor $T:W \to X(a,b)$.
Let $q: H \to \Z/3 \times \Z/3$ be the group
quotient homomorphism $q = a_{12} \times a_{23}$ and let $\kappa^{X(a,b)}: \sO_{X(a,b)}^* \to H^1(X(a,b),\Z/3)$
be the Kummer map obtained from $H^*(X(a,b),-)$.
\begin{lm}\label{lm:kappa}
  There is an isomoprhism of $\Z/3\times\Z/3$-torsors
  \[
    q_* W \cong \kappa^{X(a,b)}(a) \times \kappa^{X(a,b)}(b).
  \]
\end{lm}
\bpf
  As a $\Z/3\times\Z/3$-torsor, $\kappa^{X(a,b)}(a)\times\kappa^{X(a,b)}(b)$ comes from pulling
  back
  \[
    \Spec k[z_1,\dotsc,z_{10},\gamma,\delta]/(I,\gamma^3-a,\delta^3-b)
      \to\Spec k[z_1,\dotsc,z_{10}]/I
  \]
  by the open immersion $X(a,b)\to\Spec k[z_1,\dotsc,z_{10}]/I$. Notice that
  $W=\Spec(k[V]\otimes_k k[\alpha^{\pm1},\beta^{\pm1}]) - T\inv(S)$, so
  $q_*W=(W\times(\Z/3\times\Z/3))/H$ is an open subset of
  \[
    \Spec \bigg(\prod_{\Z/3\times\Z/3}k[V]\otimes_k k[\alpha^{\pm1},\beta^{\pm1}]\bigg)^H.
  \]

  For any $(m,n)\in\Z/3\times\Z/3$, we can map
  \[
    k[z_1,\dotsc,z_{10},\gamma,\delta]/(I,\gamma^3-a,\delta^3-b)
      \to k[V]\otimes_k k[\alpha^{\pm1},\beta^{\pm1}]
  \]
  by $\gamma\mapsto\zeta^m\alpha$, $\delta\mapsto\zeta^n\beta$, and by $\theta$
  (see \eqref{eq:theta_def}) on $k[z_1,\dotsc,z_{10}]$. This gives a map
  \[
    f:k[z_1,\dotsc,z_{10},\gamma,\delta]/(I,\gamma^3-a,\delta^3-b)
      \to\prod_{\Z/3\times\Z/3}k[V]\otimes_k k[\alpha^{\pm1},\beta^{\pm1}]
  \]
  which takes monomials
  \[
    g(z_1,\dotsc,z_{10})\gamma^i\delta^j\mapsto
      \big(\theta(g)\zeta^{im+jn}\alpha^i\beta^j\big)_{(m,n)}.
  \]

  We claim that $\im(f)$ lies in the $H$-fixed points. Suppose we have $h\in H$
  and a monomial
  \[
    g(z_1,\dotsc,z_{10})\gamma^i\delta^j\in
      k[z_1,\dotsc,z_{10},\gamma,\delta]/(I,\gamma^3-a,\delta^3-b).
  \]
  Then, using that $\im(\theta)$ is fixed under the action of $H$,
  \[
    h\cdot\big(\theta(g)\zeta^{im+jn}\alpha^i\beta^j\big)_{(m,n)}
      = \big(\theta(g)\zeta^{i(m+a_{12}(h))+j(n+a_{23}(h))}\alpha^i\beta^j\big)_{(m+a_{12}(h),n+a_{23}(h))}.
  \]
  This shows that $\im(f)$ is fixed under the action of $H$. Hence we get a map on schemes
  \[
    \Spec\bigg(\prod_{\Z/3\times\Z/3}k[V]\otimes_k k[\alpha^{\pm1},\beta^{\pm1}]\bigg)^H
      \to\Spec k[z_1,\dotsc,z_{10},\gamma,\delta]/(I,\gamma^3-a,\delta^3-b).
  \]
  This scheme map is also a map of $\Z/3\times\Z/3$-torsors, which we can
  pullback to a map $\kappa^{X(a,b)}(a)\times\kappa^{X(a,b)}(b)\to q_*W$
  of $\Z/3\times\Z/3$-torsors. Since any map of $\Z/3\times\Z/3$-torsors is an
  isomorphism, we are done.
\epf

We use these torsors to prove a theorem, but we first prove a lemma that will
be useful in proving the theorem.

\begin{lm}\label{lm:vanishing_cup}
  Let $a,b\in F^*$.
  If $\tau:\Gal F\to H$ is such that $q\circ\tau=\kappa(a)\times\kappa(b)$, then
  $\kappa(a)\smile\kappa(b)=0$ in $\rH^2(\Spec F,\Z/3)$.
\end{lm}

\bpf
  Since $q\circ\tau=\kappa(a)\times\kappa(b)$, we can write
  \begin{equation}\label{eq:G_homom}
    \tau(g) = \begin{pmatrix}
      1 & \kappa(a)(g) & \tau_{13}(g) \\
      0 & 1 & \kappa(b)(g) \\
      0 & 0 & 1
    \end{pmatrix},
  \end{equation}
  where $\tau_{13}: H \to \Z/3$ simply assigns the $(1,3)$-entry of $\tau(g)$.
  Notice that since $\tau$ is a homomorphism, we can calculate $\tau_{13}(g_1g_2)$
  for any $g_1,g_2 \in H$ -- the $(1,3)$-entry of $\tau(g_1) \tau(g_2)$ is
  \begin{equation} \label{eq:tau_13}
    \tau_{13}(g_1 g_2) = \tau_{13}(g_1) + \tau_{13}(g_2) + \kappa(a)(g_1) \kappa(b)(g_2).
  \end{equation}

  For $d:C^1(H,\Z/3) \to C^2(H,\Z/3)$ the differential map, we claim that
  $d(\tau_{13}) = - \kappa(a) \smile \kappa(b)$, whence $\kappa(a)
  \smile \kappa(b) = 0$ in $H^2(\Spec F, \Z/3)$. By \eqref{eq:tau_13}, we have
  \begin{align*}
    d(\tau_{13})(g_1,g_2) &= g_1 \tau(g_2) - \tau_{13}(g_1 g_2) + \tau_{13}(g_1) \\
      &= g_1 \tau_{13}(g_2) - \tau_{13}(g_2) - \kappa(a)(g_1)\kappa(b)(g_2).
  \end{align*}
  But the action of $\Gal F$ on $\Z/3$ is trivial; hence
  \[
    d(\tau_{13})(g_1,g_2) = - \kappa(a)(g_1) \kappa(b)(g_2)
      = - (\kappa(a) \smile \kappa(b))(g_1,g_2)
  \]
  as desired.
\epf

\begin{tm} \label{thm:main_theorem}
  For $a,b \in F^*$, the scheme $X(a,b)$ has an $F$-point if and only if $\kappa(a) \smile \kappa(b) = 0$ in $\rH^2(\Spec F,\Z/3)$.
\end{tm}

This says that $X(a,b)$ is a splitting variety for $\kappa(a) \smile \kappa(b)$.
One direction of this proof is nearly identical to \cite{HW}, and I thank
Kirsten Wickelgren for going over that proof with me.

\bpf
  Suppose there is an $F$-point $f: \Spec F \to X(a,b)$. Then
  the pullback $f^*W$ is an $H$-torsor. Using Lemma \ref{lm:kappa} and the
  fact that pushforward and pullback commute,
  \[
    q_* f^* W = f^*(\kappa^{X(a,b)}(a) \times \kappa^{X(a,b)}(b)) = \kappa(a) \times \kappa(b),
  \]
  where the last equality follows from naturality of the Kummer map.
  Viewing $f^*W$ as an element of $\rH^1(\Gal F, H)$,
  let $\tau: \Gal F \to H$ be a representative of $f^*W$. Then $q_* f^* W = \kappa(a) \times \kappa(b)$
  viewed as an element of $\rH^1(\Gal F, \Z/3\times\Z/3)$ has
  representative $q \circ \tau: \Gal F \to \Z/3 \times \Z/3$ which takes
  $g \mapsto (\kappa(a)(g),\kappa(b)(g))$, and this is the only representative.
  Thus by Lemma \ref{lm:vanishing_cup}, we find that $\kappa(a)\smile\kappa(b)=0$.

  Conversely, suppose $\kappa(a) \smile \kappa(b) = 0$. Then by \cite[Thm 2.4]{Dwy} there is an element
  of $\rH^1(\Spec F, H)$ with representative $\sigma: \Gal(F) \to H$ such that
  \[
    q_* \sigma = \kappa(a) \times \kappa(b).
  \]
  Take any cube roots $\alpha$ and $\beta$ of $a$ and $b$, respectively, in
  $\overline F$. Then for any $g \in \Gal F$,
  \[
    a_{12} \sigma (g) = \frac{g(\alpha)}{\alpha} \in \Z/3.
  \]
  Similarly, $a_{23} \sigma (g) = g(\beta)/\beta$.

  Let $\eta: H \to \GL_3(F)$ be the homomorphism of the representation $\Ind^H_N \rho$
  of $H$ defined earlier. Since $\rH^1(\Spec F, \GL_3(\overline F))$ is the pointed set
  with one element, the image of $\sigma$ under $\eta_*:\rH^1(\Spec F, H) \to
  \rH^1(\Spec F, \GL_3(\overline F))$ is trivial. Hence there is $A\in\GL_3(\overline F)$
  such that
  \[
    \eta \sigma(g) = A\inv  (gA)
  \]
  for all $g \in \Gal F$.

  Let $\pi_i: \overline{F}^3 \to \overline F$ be projection onto the $i$th
  coordinate. Consider the linear maps $A_i: F^3 \to \overline F$ which take $\mu \mapsto
  \pi_i(A\inv \mu)$. Note that $\dim \Ker(A_i) < 3$. Since $F$ is infinite, $F^3$
  is not contained in a finite union of smaller dimensional vector spaces.
  Hence there exists $\mu \in F^3$ such that $A\inv  \mu \in (\overline{F}^*)^3$. 

  Because $A\inv \mu\in(\overline F^*)^3$,
  $A\inv \mu \times (\alpha,\beta)$ determines an $\overline F$-point
  of $\G_m^3 \times \Spec k[\alpha,\alpha\inv ,\beta,\beta\inv ]$ by construction.
  For $g\in\Gal F$, notice that $g\inv (A\inv(gA)) = (g\inv A\inv )A$. Hence
  \[
    (g\inv A\inv) A = g\inv(A\inv (gA)) = g\inv(\eta\sigma(g)) = \eta\sigma(g),
  \]
  where the last step follows because the image of $\eta$ is actually contained
  in $\GL_3(F)$. Thus
  \begin{align*}
    g\inv(A\inv \mu \times (\alpha,\beta)) &= \eta\sigma(g)A\inv\mu\times (g\inv(\alpha),g\inv(\beta)) \\
      &= \eta\sigma(g)A\inv\mu\times (\alpha \kappa(a)(g\inv)(\alpha),\beta\kappa(b)(g\inv)(\beta)) \\
      &= \eta\sigma(g)A\inv\mu\times (\alpha \kappa(a)(g)(\alpha),\beta\kappa(b)(g)(\beta)) \\
      &= \eta\sigma(g)A\inv\mu\times q\sigma(g) (\alpha,\beta).
  \end{align*}
  In particular, $A\inv \mu \times (\alpha,\beta)$ determines an $\overline F$-point
  of $X(a,b)$ which by Lemma \ref{lm:fixed_ring} is fixed by the action of $\Gal F$.
  Hence $A\inv \mu \times (\alpha,\beta)$ determines an $F$ point of $X(a,b)$.
\epf

} 

\section{Automatic Realization of Galois Groups}\label{sec:galois}
The inverse Galois problem asks, for a given field $F$ and finite group $G$,
whether there is a $G$-Galois field extension $L/F$. In the affirmative case,
we say that $G$ is \textit{realizable} over $L$. A related question is that
of automatic realization: given two finite groups $G$ and $G'$, can we decide
whether $G'$ is realizable over $F$ solely from the realizability of $G$ over
$F$? An overview of the theory of automatic realizations can be found in
\cite{Jen97}.
We prove an automatic realization theorem using the results of Section
\ref{sec:splitting_variety}.

Let $G$ be a finite group. A $G$-Galois ring extension $B/A$ is then the same
thing as a $G$-torsor $\Spec B\to\Spec A$.
We prove a general lemma on the pushforward of such a torsor.

\begin{lm}\label{lm:pushforward}
  Let $B/A$ be a finite Galois extension of rings with finite Galois
  group $G$ and let $q:G\to Q$ be a quotient of $G$ with kernel $K$. Then
  $q_*\Spec B\cong\Spec A^K$.
\end{lm}

\bpf
  By the definition of pushforward of torsors,
  \[
    q_*\Spec B = \Spec\bigg(\prod_Q B\bigg)^G,
  \]
  where $G$ acts on $\prod_Q B$ by
  \[
    g\big((l_a)\big)_{a\in Q} = (gl_{q(g)\inv a})_a.
  \]
  It is enough to show that $(\prod_Q B)^G\cong B^K$. Each element of
  $Q\cong G/K$ has a representative $g\in G$. We claim that the map $B^K\to B$
  given by $x\mapsto gx$ depends only on $q(g)$. Suppose that $q(g)=q(h)$; then
  $h\inv g\in K$. But then $h\inv gx=x$ for all $x\in B^K$, so $gx=hx$. Thus
  the map induces a ring homomorphism $\phi:B^K\to\prod_Q B$ defined by
  \[
    \phi(l) = (gl)_{q(g)\in Q}.
  \]
  Note that $\phi$ is injective because $gx=0$ if and only if $x=0$. It remains
  to show that $\phi$ is surjective.

  Let $l\in B^K$. Then
  \[
    h\phi(l)=\big(hl_{q(h\inv g)}\big)_{q(g)}
  \]
  for any $h\in G$. But there is some $g'\in G$ and $k\in K$ such that
  $h\inv g=g'k$. Hence
  \[
    h\phi(l) = \big(hg'l\big)_{q(g)} = (gl)_{q(g)},
  \]
  so $\im(\phi)\subset\big(\prod_Q B\big)^G$.

  Now suppose that $(l_{q(g)})_{q(g)\in Q}$ is fixed by $G$. For any $h\in G$,
  \[
    h\inv(l_{q(g)})_{q(g)} = (h\inv l_{q(h\inv g)})_{q(g)};
  \]
  in particular, we must have that $h\inv l_{q(h)}=l_{q(1)}$. Moreover,
  $l_{q(1)}\in B^K$; for if $k\in K$, then $k(l_{q(g)})_{q(g)} = (kl_{q(g)})_{q(g)}$.
  Therefore $(l_{q(g)})_{q(g)}=\phi(l_{q(1)})$.
\epf

We now put this in the context of the rest of the paper.
Let $H$ be the mod 3 Heisenberg group and $q:H\to\Z/3\times\Z/3$ be the
quotient homomorphism. Let $F$ be a number field containing a cube root of
unity. Let $X(a,b)$ be as in Definition \ref{df:splitting_variety}.

\begin{tm}
  Suppose $a,b\in F^*$ are such that $F(\sqrt[3] a,\sqrt[3] b)/F$ is a
  $\Z/3\times\Z/3$-Galois extension. Then the following are equivalent:
  \begin{enumerate}
    \item There exists a $\Z/3$-Galois extension $L/F(\sqrt[3] a,\sqrt[3] b)$
      such that $L/F$ is an $H$-Galois extension;
    \item $\kappa(a)\smile\kappa(b)=0$ in $\rH^2(\Spec F,\Z/3)$;
    \item The scheme $X(a,b)$ has an $F$-point.
  \end{enumerate}
\end{tm}

\bpf
  (1) $\Rightarrow$ (2): Suppose that there exists a Galois $\Z/3$-extension
  $L/F(\sqrt[3] a,\sqrt[3] b)$ such that $L/F$ is a Galois $H$-extension. Then
  $\Spec L\to\Spec F$ is an $H$-torsor and
  $\Spec F(\sqrt[3] a,\sqrt[3] b)\to\Spec F$ is a $\Z/3\times\Z/3$-torsor. Let
  $\sigma:\Gal F\to H$ be a representative of $\Spec L\to\Spec F$ viewed as an
  element of $\rH^1(\Spec F,H)$. Observe
  $\kappa(a)\times\kappa(b):\Gal F\to\Z/3\times\Z/3$ is the unique cocycle
  representative of $\Spec F(\sqrt[3] a,\sqrt[3] b)\to\Spec F$. By Lemma
  \ref{lm:pushforward}, we get that
  \[
    q_*\Spec L\cong\Spec L^{\Z/3} = \Spec F(\sqrt[3] a,\sqrt[3] b);
  \]
  hence $q_*\sigma=\kappa(a)\times\kappa(b)$. Applying Lemma
  \ref{lm:vanishing_cup}, we conclude that $\kappa(a)\smile\kappa(b)=0$.

  (2) $\Rightarrow$ (3): This is Theorem \ref{thm:main_theorem}.

  (3) $\Rightarrow$ (1): Suppose there is an $F$-point $f:\Spec F\to X(a,b)$.
  Then as in the proof of Theorem \ref{thm:main_theorem}, we can pullback $f$
  to an $H$-torsor $f^*W\to\Spec F$ such that
  $q_*f^*W=\kappa(a)\times\kappa(b)$. Since $F(\sqrt[3]a,\sqrt[3]b)/F$ is a
  $\Z/3\times\Z/3$-Galois extension, it follows that
  $q_*f^*W=\Spec F(\sqrt[3]a,\sqrt[3]b)$. Note that $f^*W\to\Spec F$ is finite
  by \cite[V Prop 2.6 (iii) (i)]{SGA1}, so $f^*Y=\Spec L$ for a finitely generated
  $F$-module $L$. Furthermore, $L$ has an action of $H$ and $L^H=F$. Thus it is
  enough to show that $L$ is a field.

  For this, it suffices to show that $L$ is
  an integral domain, because $L$ is a finitely generated $F$-module. By Lemma
  \ref{lm:pushforward},
  \[
    L^{\Z/3} = q_*f^*W = F[\sqrt[3]a,\sqrt[3]b].
  \]
  Hence $\Spec L$ is a $\Z/3$-torsor over $\Spec F[\sqrt[3]a,\sqrt[3]b]$. But
  these are all known: in particular, either $L$ is a field or this is the
  trivial torsor. We show that $\Spec L\to\Spec F[\sqrt[3]a,\sqrt[3]b]$ is not
  trivial, which completes the proof.

  Suppose it were trivial, that is
  \[
    \Spec L \cong \Spec \prod_{\Z/3} F[\sqrt[3]a,\sqrt[3]b].
  \]
  We show that this implies that $\Spec L\to\Spec F$ is not an $H$-torsor,
  which gives a contradiction. We know that
  $E_{13}$ acts on $L$ by permuting the factors. We can use this to solve for
  the actions of $E_{12}$ an $E_{23}$. For $x\in F[\sqrt[3]a,\sqrt[3]b]$, the
  action of $E_{12}$ on $(x,x,x)$ is given by acting $E_{12}$ on each
  coordinate. Let $E_{12}(1,0,0)=(d,e,f)$. Then
  \[
    E_{12}(0,1,0) = E_{13} E_{12} (1,0,0) = (f,d,e)
  \]
  because the actions of $E_{23}$ and $E_{13}$ commute on $\prod_{\Z/3} F$.
  Similarly, $E_{12}(0,0,1)=(e,f,d)$.

  Let $g=d+e+f$. Then $E_{12}(1,1,1)=(g,g,g)$. But $H$ fixes $F$ and $E_{12}$
  acts on $(1,1,1)$ by acting on each coordinate, so $g=1$. Now,
  $(1,0,0)(0,1,0)=(0,0,0)$, so $E_{12}(1,0,0)(0,1,0)=(0,0,0)$. On the other
  hand,
  \[
    E_{12}(1,0,0)(0,1,0)=(d,e,f)(f,d,e)=(df,de,ef).
  \]
  Hence two of $d,e,f$ are $0$ and the other is $1$.

  From this and the fact that $E_{12}(x,0,0)=E_{12}(1,0,0)(x,x,x)$, we find
  that $E_{12}(x,0,0)$ is one of $(E_{12}x,0,0)$, $(0,E_{12}x,0)$, or
  $(0,0,E_{12}(x))$. This choice will determine the entire action
  of $E_{12}$ on $L$. The same argument works for $E_{23}$ as well.

  We have shown that $E_{12}$ and $E_{23}$ act by translating the factors and
  then acting on each factor. But then it follows that $E_{12}$ and $E_{23}$
  commute, which is a contradiction. Therefore
  $\Spec L\to\Spec F[\sqrt[3]a,\sqrt[3]b]$ is not the trivial torsor.
\epf

\bibliographystyle{amsalpha2} \bibliography{Splitting_variety_mod_3_cup_product}

\end{document}